\documentclass[12pt, amssymb, amsmath]{article}
\usepackage{epsfig,amssymb,epic,eepic,color}
\usepackage{graphicx}

\begin{document}

\renewcommand{\thesection}{\arabic{section}}
\renewcommand{\thesubsection}{\thesection.\arabic{subsection}}


\newtheorem{lem}{Lemma}[section]
\newtheorem{propo}[lem]{Proposition}
\newtheorem{theo}[lem]{Theorem}
\newtheorem{rema}[lem]{Remark}
\newtheorem{remas}[lem]{Remarks}
\newtheorem{coro}[lem]{Corollary}
\newtheorem{defin}[lem]{Definition}
\newtheorem{hypo}[lem]{Hypoth\`ese}
\newtheorem{exem}[lem]{Exemple}
\newtheorem{conj}[lem]{Conjecture}

\newcommand{\real}{{\bf R}}
\newcommand{\cn}{{\bf C}^{n}}
\newcommand{\reall}{{\bf R}^{p}}
\newcommand{\realll}{{\bf R}^{k}}
\newcommand{\mreal}{M\times{\bf R}^{p}}
\newcommand{\mreall}{M\times{\bf R}^{2k}}
\newcommand{\ri}{\rightarrow}
\newcommand{\tast}{T^{\ast}M}
\newcommand{\tastl}{T^{\ast}L}
\newcommand{\tastilde}{T^{\ast}\widetilde{M}}
\newcommand{\cd}{D_{\bullet}}
\newcommand{\cdk}{D_{\bullet-k}}
\newcommand{\cfx}{C_{\bullet}(f,\xi )}
\newcommand{\cald}{{\cal D}(M)}
\newcommand{\cm}{C_{\bullet}(M)}
\newcommand{\cmk}{C_{\bullet-k}(M)}
\newcommand{\cov}{\widetilde{M}}
\newcommand{\calm}{{\cal M}}
\newcommand{\cala}{{\cal A}}
\newcommand{\calc}{{\cal C}}
\newcommand{\call}{{\cal L}}
\newcommand{\caly}{{\cal Y}}
\newcommand{\gol}{\, \, \, \, \, \, \, \, }
\newcommand{\hatf}{\hat{f}}
\newcommand{\hatx}{\hat{x}}
\newcommand{\hatv}{\hat{v}}
\newcommand{\mn}{M^{n}}
\newcommand{\barm}{\bar{M}}
\newcommand{\bfs}{{\bf S}}
\newcommand{\piu}{\pi_{1}}
\newcommand{\pii}{\pi_{i}(M)}
\newcommand{\omi}{\Omega(L_{0},L_{1})}
\newcommand{\omibar}{\Omega(\bar{L}_{0},\bar{L}_{1})}
\newcommand{\omil}{\Omega(L,L)}
\newcommand{\omibarl}{\Omega(\bar{L},\bar{L})}
\newcommand{\la}{^}
\newcommand{\bfc}{{\bf C}}
\newcommand{\covl}{\widetilde{L}}
\newcommand{\barl}{\bar{L}}
\newcommand{\tdn}{{\bf T}^{2n+1}}
\newcommand{\tn}{{\bf T}^{n}}
\newcommand{\bft}{{\bf T}}
\newcommand{\bfz}{{\bf Z}}

\begin{center}
{\large\bf  ON THE TOPOLOGY OF MONOTONE LAGRANGIAN SUBMANIFOLDS}\\
\vspace{.3in}
\noindent Mihai DAMIAN\\
\noindent Universit\'e de  Strasbourg\\
IRMA, 7, rue Ren\'e Descartes,\\
67 084 STRASBOURG\\
e-mail : damian@math.u-strasbg.fr\\

\end{center}
\vspace{.4in}

\noindent{\bf Abstract\, :}\, We find new obstructions on the topology of monotone Lagrangian submanifolds of $\cn$ under some hypothesis on the homology of their universal cover. In particular we show that nontrivial connected sums of manifolds of odd dimensions do not admit monotone Lagrangian embeddings into $\cn$ whereas some of these examples are known to admit usual  Lagrangian embeddings: the question of the existence of a monotone embedding for a given Lagrangian in $\cn$ was open. In dimension three we get as a corollary that the only orientable Lagrangians in ${\bf C}^{3}$ are products $\bfs^{1}\times \Sigma$. The main ingredient of our proofs is the lifted Floer  homology theory which we developed in \cite{Da2}. 
\vspace{.2in}

\section{Introduction}

\subsection{Motivation}

The starting point of this paper is the following general question :

\begin{quote} What can be said about the topology of a closed Lagrangian $L \subset \cn $? 
\end{quote}

A celebrated theorem of M. Gromov \cite{G} asserts that such a submanifold cannot be simply connected and it is easy to prove that it also  has to satisfy $\chi(L)=0$ if it is orientable (this equality only holds modulo $4$ in general \cite{Au}), but can one say more?  The answer is affirmative  in low dimensions;  for instance when $n=2$  we know precisely the surfaces admitting Lagrangian embeddings: the only orientable one is the torus, and the non-orientable Lagrangians  are those whose  Euler characteristic is a multiple of $4$ \cite{Giv}, with the exception of the Klein bottle \cite{Sch}, \cite{N}. We also know pretty much about the closed  Lagrangians in ${\bf C}^{3}$: a quite recent theorem of K. Fukaya \cite {Fu} (2006) asserts:

\begin{theo}\label{Fu} If $L$ is a three-manifold which is closed orientable and prime (i.e. not decomposable into a nontrivial connected sum) then it admitts a Lagrangian embedding into ${\bf C}^{3}$ if and only if it is diffeomorphic to a product $\bfs^{1}\times \Sigma$. \end{theo}

In greater dimensions it is of course too optimistic to expect such a precise characterisation of the topology of Lagrangian submanifolds. There are not many examples of Lagrangians available and we do not know if this lack of examples comes from some strong topological restrictions still to be discovered. One of the first papers concerned with the topology of closed Lagrangians in $\cn$ and the possible examples was written by M. Audin \cite{Au} in 1988. The author observed that all the known orientable Lagrangian submanifolds have a common feature: they fiber over the circle. This led to the natural question of the existence of other examples.  Few years later, in 1991, L. Polterovich \cite{Pol}  gave a negative answer to this question: he constructed a lot of new examples of Lagrangian submanifolds starting from Lagrangian immersions and replacing neighborhoods of the double points by $1$-handles. The resulting Lagrangians are often connected sums, which do not fiber over the circle. Indeed, it can be proved  more generally that a manifold whose fundamental group is a nontrivial free product, cannot fiber over the circle (see \cite{Da1}, Prop.2.3 for instance). Here is the theorem of Polterovich \cite{Pol}: 

\begin{theo}\label{Pol} Let $P=\bfs^{n-1}\times \bfs^{1}$ and $Q$ the manifold obtained from $\bfs^{n-1}\times [0,1]$ after gluing the points $(x,1)$ and $(\tau x, -1)$ where $\tau$ is the standard reversing orientation involution on $\bfs^{n-1}$. Then:\\
a) Let $L_{1}$ and $L_{2}$ be closed connected manifolds admitting Lagrangian embeddings into $\cn$. Then $L_{1}\#L_{2}\#Q$ admits a Lagrangian embedding into $\cn$. Moreover $L_{1}\#L_{2}\#P$ admits a Lagrangian embedding into $\cn$ if $n$ is odd.\\
b) Let $L$ be a closed connected manifold admitting a Larangian immersion into $\cn$. Then $L\#kQ$ admits a Lagrangian embeding into $\cn$ for some non-negative integer $k$. Moreover $L\#kP$ admits a Lagrangian embedding into $\cn$ if $n$ is odd. 
\end{theo}

L. Polterovich also points out  that any closed orientable $3$-manifold admits a Lagrangian immersion in $\bfc^{3}$. This  gives many orientable Lagrangians  (of odd dimension) not fibering over the circle, as $\tdn\#\tdn \#P$, $P\#P\# P$, $\tdn\#P\#P$ etc.  Taking cartesian products, we get examples of even dimension too.  \\

The developement of Floer homology techniques led to an intensive study of a special class of Lagrangian submanifolds called monotone. Recall that  a Lagrangian submanifold $L\subset (M,\omega)$ is defined to be  monotone if the morphism $I_{\omega}:\pi_{2}(M,L)\ri \real$ given by the symplectic area and the morphism $I_{\mu}:\pi_{2}(M,L)\ri {\bf Z}$ defined by the Maslov index \cite{Ar} are positively proportional. Monotone submanifolds are known to be more rigid with respect to Lagrangian intersections. Also some constraints on their topology were established (\cite{Bi}, \cite{BiC}, \cite{BiCo}, \cite{BiCo2}, \cite{BiCo3}) but     these topological properties are not  known to be specific to monotone Lagrangians. It is therefore very tempting to reformulate Audin's question in these terms: 

\begin{quote} {\bf Question 1.} Let $L\subset \cn$ be a closed orientable monotone Lagrangian. Does $L$ fiber over $\bfs^{1}$?
\end{quote}

In view of Polterovich's examples another natural question related to the previous is the following: 

\begin{quote} {\bf Question 2.} Let $L\subset \cn$ a closed (orientable) Lagrangian. Does $L$ also admit a monotone Lagrangian embedding into $\cn$? 
\end{quote} 

The aim of this paper is to show that there is indeed a topological rigidity specific to monotone Lagrangian submanifolds. The answer to Question~2 is negative, and many of  Polterovich's examples cannot be embedded as monotone Lagrangians. We are able to prove a weaker version of the assertion of Question 1. (``stable fibration'' instead of fibration), but only under a hypothesis related to the holomorphic disks  of Maslov index equal to $2$, with boundary in $L$ (in particular the Maslov number of these submanifolds is supposed to be $N_{L}=2$). Note that a topological constraint for monotone Lagrangians with $N_{L}=n$ was established in our previous paper (\cite{Da2}, Th.1.7).  Recall that the Maslov number $N_{L}\in {\bf N}$ is defined to be the positive generator of $\mbox{Im}(I_{\mu})$.

The case of $3$-Lagrangians is easier to deal with. In this case we prove that the answer to Question 1 is positive; moreover we have $L=\bfs^{1}\times \Sigma$, as in Fukaya's theorem. As far as we know  this question is completely open in greater dimensions. 

\subsection{Main results} 

All our results are valid for closed orientable Lagrangians $L\subset M$  in arbitrary monotone symplectic manifolds, with the property that $\phi(L)\cap L=\emptyset$, for some Hamiltonian diffeomorphism $\phi$: these submanifolds are called displaceable. Of course any closed submanifold in $\cn$ or more generally in $M=\bfc\times W$,  is displaceable. We suppose without restricting generality that all the Lagrangian submanifolds we consider are of dimension greater than three.  Here is our first result: 

\begin{theo}\label{general} Let $L\subset M$ be a closed orientable Lagrangian which is displaceable. Denote by $\widetilde{L}$ its universal cover. Suppose that 

(a) $H_{2i+1}(\widetilde{L}, {\bf Z}/2)=0$ for all $i$.

(b) $\piu(L)= G_{1}\ast G_{2}$ with $\mbox{rk}(G_{j}/[G_{j},G_{j}])\neq 0$ for $j= 1, \, 2$.\\
Then $L$ is not monotone. 
\end{theo}

Here $\ast$ denotes the free product of two groups. As one can easily check, in odd dimension the hypothesis (a) is preserved by connected sums. We have: 

\begin{theo}\label{connected} Let $L^{n}\subset M^{2n} $ a closed orientable Lagrangian of odd dimension which is displaceable. Suppose that $L=L_{1}\#L_{2}$ and the manifolds $L_{j}$ satisfy the hypothesis (a) of \ref{general}. Then $L$ is not monotone.

\end{theo}

We get thus many examples of closed manifolds of odd dimension which do not admit monotone Lagrangian embeddings into $\cn$: arbitrary connected sums of manifolds of the form  $K(\pi, 1)$, $ \bfs^{2k}\times\bfs^{1}$, or any other summand satisfying the hypothesis (a). In particular a lot of Polterovich's examples   -- the odd-dimensional orientable ones obtained from Th.\ref{Pol} -- do not admit monotone embeddings. 

Also note that hypothesis (a) is satisfied by any $3$-manifold with infinite fundamental group. Applied to $3$-manifolds, the previous theorem implies that monotone orientable $3$-manifolds must be necessarily prime. In fact we have a more precise result: 

\begin{theo}\label{product}  Let $L^{3}\subset \bfc^{3}$ be a closed orientable monotone Lagrangian, then $L$ is diffeomorphic to $\bfs^{1}\times \Sigma$. 
  
\end{theo}

This theorem can be seen as a corollary of Th.\ref{general} and Fukaya's theorem \ref{Fu}. Actually we are able to give another proof which avoids the technical difficulties of Fukaya's proof. Another proof of  Th.\ref{product}, was independently obtained using a different method by J. D.  Evans and J. Kedra  \cite{EK}. \\

As a matter of fact, all the results we stated in this section are consequences of a more general one which we will present below. 

In order to  do this  we need to introduce the group ring 
$\Lambda={\bf Z}/2\, [\piu(L)]$ of the fundamental group of a manifold, as well as its completion, the Novikov ring $\Lambda_{u}$ associated to some morphism $u:\piu(L)\ri\real$:
$$\Lambda_{u}=\{\sum_{i\geq 1}n_{i}g_{i}, \, n_{i}\in {\bf Z}/2, \, g_{i}\in \piu(L), \, u(g_{i})\ri\infty\}.$$
By definition the Novikov homology $H_{\ast}(L,u)$ is given by 
$$H_{\ast}(L;u)\, =\, H_{\ast}(\Lambda_{u}\otimes_{\Lambda} C_{\bullet}^{\widetilde{L}}(f,\xi)),$$
where $f$ is a Morse function on $L$, $\xi$ is a generic gradient, and $C_{\bullet}^{\widetilde{L}}$ is the associated ($\Lambda$-free) Morse complex obtained by lifting the gradient lines to the universal cover. We will prove:
\begin{theo}\label{moregeneral} Let $L\subset M$ be a closed monotone orientable Lagrangian which is displaceable. Suppose that 

(a)  $H_{2i+1}(\widetilde{L}, {\bf Z}/2)=0$ for all $i$.

Then we have :

(i)  The Novikov homology $H_{\ast}(L;u)$ is  $\Lambda$-torsion for any morphism $u:\piu(L)\ri\real$. 

(ii) There is some element $g\in \piu(L)$ which has a finite number of conjugacy classes. In particular $\piu(L)$ admits a finite index subgroup with non-trivial center. 

\end{theo}

 In even dimensions the (known) fact that $\chi(L)~=~0$ for any displaceable Lagrangian $L$ is a stronger restriction on the topology of connected sums (since Euler charactestic is not additive). We are able to adapt the statement of \ref{connected} in this case, but we cannot apply it to a concrete example. Instead, we will use \ref{moregeneral} to  prove that cartesian products of Polterovich's examples are  counter-examples for Question 2. 
Note that condition (a)  is also preserved by cartesian products. We have: 

\begin{coro}\label{pair} Let $L_{1}$, $L_{2}$ be two closed odd-dimensional orientable manifolds which are both nontrivial  connected sums of  manifolds satisfying (a) of Th.\ref{general}. Then $L=L_{1}\times L_{2}$ does not admit any monotone displaceable Lagrangian  embedding.
\end{coro}

As a consequence, arbitrary products of the orientable examples produced by Th.\ref{Pol}, do not admit monotone embeddings into $\cn$ either, whereas they obviously admit Lagrangian embeddings. For instance 
$$ L=[(\bfs^{2n}\times \bfs^{1})\#\tdn\#\tdn]\times [(\bfs^{2n}\times\bfs^{1})\#\tdn\#\tdn]$$ does not admit any monotone embedding into ${\bf C}^{8n+4}$. \\

Let us now state a result which is related to Question 1. Consider a closed monotone Lagrangian $L^{n}\subset (M^{2n}, \omega)$ whose Maslov number is $N_{L}=2$. Let $J$ be an almost complex structure on $M$ which is compatible with $\omega$. Let $[g]\in H_{1}(L, \bfz)$ and define  the following space of $J$-holomorphic disks with boundary in $L$ : 
$$\calm_{[g]}(M,L,J,2)\,  = \, \{\, w:(D, \partial D)\ri (M,L)\, |\, \bar{\partial}_{J}w=0, \, I_{\mu}(w)=2, \, [\partial w]=[g]\, \}.$$
Then consider 
$${\cal N}_{[g]}= \calm_{[g]}\times \bfs^{1}/\mbox {PSL}_{2}(\real), $$
where $\mbox {PSL}_{2}(\real)$ acts on $\calm_{[g]}\times \bfs^{1}$ by $$\psi\cdot (w,z)\, =\, ( w\circ \psi, \psi^{-1}(z)).$$
 
For a generic choice of $J$, both $\calm_{[g]}$ and ${\cal N}_{[g]}$ are closed manifolds of dimensions $n+2$ resp. $n$ (see \cite{Da2}). There is a natural evaluation map 
$$\mbox{ev}_{[g]}: {\cal N}_{[g]}\ri L$$ defined by $$\mbox{ev}_{[g]}([w,z])\, =\, w(z).$$
The elements of a fiber $\mbox{ev}_{[g]}^{-1}(x)$  are in one-to-one correspondence with the (unparametrized) $J$-holomorphic disks with boundary in $L$ of class $[g]$,   with Maslov number equal to two, and passing through $x$. The parity of this number does not depend on a generic choice of $x$; it can be shown that it does not depend on the choice of a generic almost complex structure either.  Let us denote by $\#_{2}^{[g]}$ the degree modulo $2$ of this evaluation map (one could also call it relative Gromov invariant). The union $\calm=\bigcup_{[g]\in H_{1}(L,\bfz)}\calm_{[g]}$ is still a closed manifold of dimension $n+2$ and we can define analogously $\cal N$ and consider the degree modulo $2$ of the evaluation map $\#_{2}$. Finally, for $g\in \piu(L,x)$ we can consider the parity $\#_{2}^{g}(x)$  of the number of holomorphic disks of Maslov number $2$ passing through  a generic $x$ with boundary in $L$ whose homotopy class is  $g$. We obviously  have:
$$\#_{2}^{[g]}\, =\, \sum_{h\in \piu(L,x), \, [g]=[h]}\#_{2}^{h}(x).$$

A generic path between two points $x$ and $y$ yields  a cobordism between $\mbox{ev}^{-1}(x)$ and $\mbox{ev}^{-1}(y)$. Notice that if two disks in $\mbox{ev}^{-1}(x)$ whose boundaries define the same element in $\piu(L,x)$ are cobordant respectively   to two disks in $\mbox{ev}^{-1}(y)$ then the homotopy classes of their boundaries coincide in $\piu(L,y)$. More precisely,    using the identification between $\piu(L,x)$ and $\piu(L,y)$ given by the same path, there is a cobordism between the Maslov $2$-disks passing through $x$ in the class $g$ and the Maslov $2$-disks passing through $y$ in the same class $g$. Therefore, when such a path between $x\neq y$  is given we may write $\#_{2}^{g}(x)=\#_{2}^ {g}(y)$. In this case we will denote the common value by $\#_{2}^{g}$. 

Let us now state our theorem. 

\begin{theo}\label{fibration} Let $(M,\omega)$ be an exact symplectic manifold and $L\subset M$ be a closed  monotone Lagrangian (not necessarily displaceable) with $N_{L}=2$. Suppose that $\#_{2}\neq 0$. Then, there is a class $u\in H^{1}(L, \bfz)$ such that $H_{\ast}(L;u)=0$.\end{theo}

If $L$ is spin then the spaces $\calm$ and $\cal{N}$ are orientable \cite{FO3}. Therefore we can define analogously $\#, \, \#^{[g]},  \, \in \bfz$, as the degrees of the respective evaluation maps, as well as $\#_{2}^{g}$.  We prove:

\begin{coro}\label{spin}
 Let $L$ be as in the previous theorem and spin. Then if $\#~=~1$, we have

(i) If   $\piu(L)$ is polycyclic then $L  \times \bfs^{3}$ (and more generally $L  \times \bfs^{2k+1}$)  fibers over $\bfs^{1}$. 

(ii) More generally,    $L  \times \bfs^{3}\times N$ fibers over $\bfs^{1}$ for any closed manifold $N$ with $\beta_{1}(N)\neq 0$. 
\end{coro}

\noindent{\bf Remark}\\
 If $L\subset M$ is a closed monotone orientable Lagrangian which is displaceable and satisfies the condition (a) of Th.~\ref {general}, then it was shown in \cite{Da2} (Theorem~\ref{dadoi} below) that $N_{L}=2$ and $\#_{2}^{g}\neq 0$ for some $g\in \piu(L)$. But this has no direct consequence on the vanishing of $\#_{2}$. 

\subsection{Idea of the proofs}

In \cite{Da2} we developed a new version of Floer homology theory for monotone Lagrangian submanifolds which involves some (arbitrary) covering space of the given Lagrangian. The outcome of this theory when $L\subset M$ is monotone, displaceable  with $N_{L} \geq 3$ is a spectral sequence whose initial page is built with the homology of the considered covering space and whose limit is zero. 

Unlike the usual Floer homology, when $N_{L}=2$, this ``lifted'' Floer homology is  not defined in general. This enables us to argue by contradiction: for instance in the proof of Th.1.5 in \cite{Da2} (which is Th.\ref{dadoi} below), if  $L$ is monotone, orientable,  displaceable and satisfies (a) of Th.\ref{general}, the supposition that lifted Floer homology is well defined leads to a contradiction and therefore to the relation $N_{L}=2$. In the present paper we look more closely to the situation where this homology {\it is not well defined} and come to the conclusion that in this case the Novikov homology $H_{\ast}(L;u)$ must be $\Lambda$-torsion for any cohomology class $u$ (compare to the statement of Th.\ref{moregeneral}).  This is also  valid for the class $u=0$ for which the conclusion is that the usual  homology of the universal cover  $H_{\ast}(\widetilde{L}, \bfz/2)$ is $\Lambda$-torsion. Moreover, we show that in this case there is a finite index subgroup of $\piu(L)$ with trivial center. 

Now we know from \cite{Da2} some situations where the lifted Floer homology is not well defined -- at this point we make use of the condition (a) which appears in the statement of Th.\ref{general}. On the other hand, we show that under some hypothesis -- as (b) in Th.\ref{general} for instance -- some of the  Novikov homology groups cannot be $\Lambda$-torsion.  Alternatively, hypothesis (b) implies that there is no finite index subgroup of the fundamental group having non-trivial center.  When both conditions (a) and (b)  are fulfilled, each of these arguments leads us to a contradiction, whose only possible explanation is that a monotone displaceable  Lagrangian embedding of the considered manifold $L$ does not actually exist. \\

The paper is organized as follows. In the next section we recall the definition of the lifted Floer homology from \cite{Da2} and analyse the case $N_{L}=2$ where it is not well defined. The last section contains the proof of results 1.3$\ldots$1.9 which are stated above. 

\section{Lifted Floer homology}

\subsection{Definition} 

Start with a monotone Lagrangian $L\subset M$ with $N_{L}\geq 3$ and  a covering space $p:\bar{L}\ri L$. Let $(\phi_{t})_{t\in [0,1]}$ be a Hamiltonian isotopy  of $M$ such that $L$ and $\phi_{1}(L)$ are transverse. Build the usual Floer complex associated to $L$ and $\phi_{1}(L)$, spanned by the intersection points $L \cap \phi_{1}(L)$. Its differential is defined by counting isolated $J$-holomorphic strips (for a fixed, generic almost complex structure $J$). Each such strip $u:\real\times [0,1]\ri M$ defines a path $\gamma_{u}(s) =u(s,0)$ on $L$. By lifting these paths to $\bar{L}$ one can define another complex, which is spanned by all the points of the fibers $p^{-1}(x)$ for $x\in L\cap \phi_{1}(L)$. Its differential is 
$$\delta \bar{x}\, =\, \sum_{\bar{y}\in p^{-1}(y)}n(\bar{x},\bar{y})\bar{y},$$
where $n(\bar{x},\bar{y})\bar{y}$ is the number (modulo $2$) of isolated holomorphic strips $u$ between $x$ and $y$, whose associated paths $\gamma_{u}$ admit  lifts $\tilde{\gamma}_{u}$ between $\bar{x}$ and $\bar{y}$. 

As in the usual case, the homology of this complex does not depend on the choices of $(\phi_{t})$ and $J$ and it is related to the singular homology of $\bar{L}$ by a spectral sequence analogous to the one established by Y-G. Oh in \cite{Oh} (see also \cite{Bi}).  The following theorem of \cite{Da2} (Th. 2.1) summarizes the properties of our lifted complex: 

\begin{theo}\label{lift} Let $L\subset M$  be monotone Lagrangian with $N_{L}\geq 3$ and $p:\bar{L}\ri L$ a covering space.  Let $(\phi_{t})$ be a Hamiltonian isotopy such that $L$ and $\phi_{1}(L)$ intersect transversally. 
Then there exist a free ${\bf Z}/2$-complex $C_{\bullet}$ spanned by
$\bigcup_{x\in L\cap\phi_{1}(L)}p^{-1}(x)$ such that :

(i)  The homology $HF(\barl)=H_{\ast}(C_{\bullet})$ does not depend on $(\phi_{t})$.

(ii) There are applications $\delta_{1},
\delta_{2}, \ldots, \delta_{\ell}, \ldots$ defined as follows and  with the following properties : 
\begin{itemize}
\item $\delta_{1}: H_{\ast}(\barl,{\bf Z}/2)\ri  H_{\ast-1+N_{L}}(\barl,{\bf Z}/2)$, satisfies $\delta_{1}\circ \delta_{1}=0.$ Denote by $H^{(1)}_{\ast}(\barl)$ the homology groups $H(H_{\ast}(\barl,{\bf Z}/2), \delta_{1})$ (with the initial Morse grading). 
\item For $\ell\geq 2$ $\delta_{\ell}: H_{\ast}^{(\ell-1)}(\barl)\ri  H_{\ast-1+\ell N_{L}}^{(\ell-1)}(\barl)$ satisfies $\delta_{\ell}\circ\delta_{\ell}=0$. The resulting homology is denoted by $HF_{\ast}^{(\ell)}(\barl)$. 
\item If $\delta_{\ell}=0$ for any $\ell\geq 1$ then 
$$H_{\ast}(C_{\bullet})\approx H_{\ast}(\barl,{\bf Z}/2).$$
\item  $HF(\barl)=0$ if and only if $HF^{(\ell)}(\barl)$ vanishes for some $\ell$. 
\end{itemize}

If $\barl=\covl$ is the universal cover of $L$ then $C_{\bullet}$ can be viewed as a free,
finite-dimensional complex over  $\Lambda = {\bf Z}/2\, [\piu(L)]$, spanned by $\phi_{1}(L)\cap L$. Let $u\in H^{1}(L, \real)$. The statements above are valid when one replaces $C_{\bullet}$ by $\Lambda_{u}\otimes_{\Lambda} C_{\bullet}$ and $H_{\ast}(\covl)$ with the  Novikov homology $H_{\ast}(L;u)$ (with $\bfz/2$ coefficients).
\end{theo}

\subsection{The case $N_{L}=2$}

As we said above, in this case the lifted Floer homology is not defined in general. This can be easily  seen  by considering $L=\bfs^{1}\subset \bfc$ and the lifted Floer homology of the univeral cover $\widetilde{L}$. As explained in \cite{Da2}, we still can define the application $\delta:C_{\bullet}\ri C_{\bullet}$ but this is not a differential in general. Let us recall the reason. Consider first the usual Floer setup. For any  generator $x\in L\cap \phi_{1}(L)$, we have 
$$\delta^{2}x\, =\, \sum_{y\in L\cap \phi_{1}(L)}m(x,y)y,$$
where $m(x,y)$ represents the number (modulo $2$) of broken $J$-holomorphic strips $(v_{1},v_{2})$ joining $x$ to $y$.  The standard method in Morse theory to show that the number of such broken strips is even is to prove that they are boundary points of an one-dimensional compact manifold $\bar\call^{1}(x,y)$ whose interior consists of unbroken strips from $x$ to $y$. Therefore the isolated broken strips arise in pairs.

And indeed this is what happens for $N_{L}\geq 3$, the associated paths $\gamma$  of the isolated broken holomorphic strips together with the paths corresponding to the unbroken strips of $\bar{\call}^{1}$ draw pictures of this form on $L$. These (contractible) pictures lift to pictures of the same form in any covering space $\barl$, which proves that we have a lifted Floer complex. \\

When $N_{L}=2$ the situation changes because  $J$-holomorphic disks of Maslov index $2$ with boundary in $L$ or $\phi_{1}(L)$ may  obstruct the compactness of $\bar{\call}^{1}$. Actually, for dimensional reasons, this can only  happen when  $x=y$. And when it happens, it may yield some "altered" Morse pictures where we have on one side a broken strip from $x$ to $x$ via $y$, on the other side a $J$-holomorphic disk passing through $x$ of Maslov index two and boundary contained in $L$ or $\phi_{1}(L)$, and between them an one-parameter family of unbroken strips from $x$ to $x$. This phenomenon can be easily visualized in the case of $L=\bfs^{1}\subset \bfc$ (see \cite{Da2}). 

  So a priori $\delta^{2}(x)\neq 0$:  it counts isolated broken strips from $x$ to $x$, so it equals the number modulo $2$  of   holomorphic disks passing through $x$, as above. However, Y-G. Oh shows that the total number of these disks is always even \cite{Oh}, which enables one to define the usual Floer homology for $N_{L}=2$.  In the case of lifted Floer homology this is no longer true. Consider the case of the universal cover $\widetilde{L}$. A broken holomorphic strip from $x$ to $x$ draws a loop on $L$, which  lifts to $\widetilde{L}$ as a path from some $\tilde{x}\in p^{-1}(x)$ to $g\tilde{x}$, where $g\in \piu(L)$ is the homotopy class of this loop. If this holomorphic strip is a boundary point of an altered Morse picture, then either the holomorphic disc on the other side has its boundary in $L$ and $g$ is the class of this boundary, or it has its boundary on $\phi_{1}(L)$ and then the loop drawn by the broken strip is nullhomotopic, i.e. $g=1$. Since the parity of the numbers of the disks of each type is the same, we infer that  
$$ (\delta^{\widetilde{L}})^{2} x\, =\, \sum_{g\in \piu(L)}\#_{2}^{g}(x)\, (1+g)x. $$ 
 Recall that $\#_{2}^{g}(x)$   counts modulo $2$  holomorphic disks of Maslov index $2$ and boundary in $L$ of class $g$ passing through a generic point $x$ (the intersection points $L\cap\phi_{1}(L)$ are chosen to be generic with respect to this map).  Now fix a point $\tilde{x}$ in each fiber $p^{-1}(x)$ of an intersection point $x\in L\cap \phi_{1}(L)$ (we have to do it anyway in order to construct the lifted Floer homology $HF(\covl)$ over $\Lambda$).  We thus have  homotopy classes for paths joining intersection points such that any loop resulting form the concatenation of such paths is trivial. As above, $\#_{2}^{g}(x)$ does not depend on $x$ and we write: 
$$(1)\gol  (\delta^{\widetilde{L}})^{2} x\, =\, \sum_{g\in \piu(L)}\#_{2}^{g}\, (1+g)x. $$
 For an arbitrary covering $\bar{L}$ with deck transformation group $G$, we should replace $(1+g)$ by $(1+\hat{g})$ in the formula above: here $\hat{g}$ is the image of $g$ through the projection $\pi:\piu(L)\ri G$.  We are now ready to prove the following 
\begin{propo}\label{torsion} Let $L\subset (M, \omega)$ be a closed monotone Lagrangian with $N_{L}=2$. Then at least one of the following assertions is true: 

(i) For any covering $\barl\ri L$ the lifted Floer homology in well defined and Th.\ref{lift} is valid. 

(ii) For any $u\in H^{1}(L, \real)$, $H_{\ast}(L;u)$ is $\Lambda$-torsion: there exists $\lambda\in \Lambda$, $\lambda\neq 0$ such that $\lambda h = 0$ for each $h\in H_{\ast}(L;u)$.
\end{propo} 

\noindent\underline{Proof}

When $(\delta^{\widetilde{L}})^{2}=0$ in the formula $(1)$ the number  of Maslov $2$-disks passing through an intersection point  $x$ in any fixed homotopy class $g$ is even. But the cobordism above shows that if this property is true for some $x$, then it is valid for {\it any} generic $y\in L$. Therefore $(\delta^{\widetilde{L}})^{2}=0$ for any Hamiltonian isotopy $\phi_{t}$:  the lifted homology of the universal cover  is well defined and the proof of \ref{lift} is analogous to the one in \cite{Da2}. If this is the case, then the formula above shows that the lifted Floer homology associated to an arbitrary covering $\barl$ is well defined as well. 

Otherwise let us consider a particular Hamiltonian isotopy $(\phi_{t})$:
 we choose it such that in a Weinstein neighborhood of $L$ we have  $\phi_{t}(L)=t\, \mbox{graph}(df)$ for some Morse function $f:L\ri \real$. Therefore   the usual Floer complex is spanned by the critical points of $f$ and graded by their Morse indices. According to \cite{Oh}, if $f$ is small enough, we have  
$$\delta=\delta_{0}+\delta_{1}+\ldots+\delta_{k}+\ldots,$$
where $\delta_{0}$ is the Morse differential defined by some generic gradient of $f$, and 
$$\delta_{k}:C_{i}\ri C_{i-1+kN_{L}}.$$
The same is obviously true for the differential of the lifted complex; when  $(\delta^{\barl})^{2}$ is zero, the application $\delta_{k}$ is used to define the differential on the $k$-th page of the spectral sequence and  in the case when this sequence is degenerate it coincides with the one in Th.\ref{lift}. When $(\delta^{\covl})^{2}$ does not vanish, the formula~$(1)$ above shows that it preserves the degree. Using that $N_{L}=2$, this implies that 
$$(\delta^{\widetilde{L}})^{2}\, =\, \delta_{0}\delta_{1}+\delta_{1}\delta_{0}.$$
Denote $$\lambda= \sum_{g\in \piu(L)}\#_{2}^{g}\, (1+g).$$
Viewing $\delta_{1}$ as a homotopy, it follows that the map $h\mapsto\lambda h$ from the lifted Morse complex $C_{\bullet}(L,f)$ to itself is nullhomotopic. Therefore it vanishes at the level of the homology $H_{\ast}(\covl, \bfz/2)$. The same is true for the lifted Novikov complex $\Lambda_{u}\otimes_{\Lambda}C_{\bullet}(\covl, f)$ and therefore the multiplication by $\lambda$ induces zero in Novikov homology. 

\hfill $ \Box$

\section{Applications}
In this section we give the proof of our results which we stated in \S1.2. Before that let us remind Th.1.5 of \cite{Da2} which will be used several times in the proofs: 
\begin{theo}\label{dadoi} Let $L\subset M$ be an orientable monotone Lagrangian which is displaceable and satisfies the condition $(a)$	of Th.\ref{general}. Then $N_{L}=2$ and the lifted Floer homology $HF(\covl)$ is not defined: for any generic $x\in L$ there is some nontrivial $g\in \piu(L)$ with $\#_{2}^{g}(x)=1$. In particular for any compatible almost complex structure $J$ through any $x\in L$ passes at least a $J$-holomorphic disk with boundary in $L$ and Maslov index $2$.
\end{theo} 

The proof follows easily from the previous section. Since $L$ is orientable $N_{L}$ must be even. If $HF(\covl)$ is well defined, then, using the notations from Th.\ref{lift}, we see that $H_{0}^{(\ell)}(\covl)=\bfz/2$ for all $\ell$. Therefore, according to the mentioned theorem, $HF(\covl)\neq 0$, contradicting the displaceability of $L$. So $HF(\covl)$ is not well defined which implies $N_{L}= 2$ and, using the formula $(1)$,  $\#_{2}^{g}(x)=1$ for some nontrivial $g$. The conclusion about the existence of the holomorphic disks also follows for a  non-generic $J$ using Gromov compacntess, as in \cite{Da2}. 

Now we are able to give the proofs:\\

\noindent\underline{Proof of \ref{general}}

Let $L$ be Lagrangian monotone  as in the statement of the theorem. Then according to Th.\ref{dadoi}, the condition $(a)$ implies that $N_{L}=2$ and the lifted Floer homology $HF(\covl)$ is not well defined (there is some $g\in \piu(L)$ such that $\#_{2}^{g}=1$). We are able to give two alternative proofs: \\
\\
\underline{First proof}

 By \ref{torsion} we have that the Novikov homology of $L$ is $\Lambda$-torsion. We will contradict the existence of $L$ by proving :
\begin{lem}\label{Novikov}
Suppose that $\piu(L)=G_{1}\ast G_{2}$ with $\mbox{rk}(G_{i}/[G_{i},G_{i}])\neq 0$. Then for some class $u\in H^{1}(L,\bfz)$ the Novikov homology group $H_{1}(L;u)$ is not $\Lambda$-torsion. 
\end{lem}

\noindent\underline{Proof of the Lemma}

The proof is quite similar to the one of Prop.2.3 in \cite{Da4} which shows that $H_{1}(L;u)\neq 0$ for any $u\neq 0$. Let $u_{i}\in Hom(G_{i}, \bfz)$, $u_{i}\neq 0$, and $u:\piu(L)\ri \bfz$ canonically defined by $u_{1}$ and $u_{2}$. Choose a presentation of $G_{i}$ with $p_{i}$ generators and $q_{i}$ relators. Then, as in the proof of Prop.2.3 in \cite{Da4} the group $H_{1}(L;u)$ can be computed from a sequence of the form: 
$$\ri\Lambda_{u}^{q_{1}}\oplus \Lambda_{u}^{q_{2}}\stackrel{\left(\begin{array}{cc}
\delta_{2}^{(1)}&0\\
0&\delta_{2}^{(2)}\end{array}\right)}{\longrightarrow}\Lambda_{u}^{p_{1}}\oplus\Lambda_{u}^{p_{2}}
\stackrel{\left(\begin{array}{cc}
\delta_{1}^{(1)}&\delta_{1}^{(2)}\end{array}\right)}{\longrightarrow}\Lambda_{u}\ri 0,$$
where on the generators $x^{(i)}$ of the chosen presentation of $G_{i}$, $\delta_{1}^{(i)}$ is defined by 
$$x^{(i)}\mapsto 1-x^{(i)}.$$
Now pick a generator $x$ of $G_{1}$ and a generator $y$ of $G_{2}$ from the given presentations. Without any  loss of generality we may suppose that $u(x)>0$ and $u(y)>0$. Denote $\mu=\sum_{j\geq 0}y^{j}=(1-y)^{-1}$. It is easy to check that $(1, -\mu(1-x))$ is in the kernel of $(\delta_{1}^{(1)}\,  \delta_{2}^{(1)})$ (to be more precise we could write this element $e_{x}-\mu(1-x)e_{y}$, where $e_{x}$ and $e_{y}$ are the elements of the canonical basis of $\Lambda_{u}^{p_{1}}\oplus \Lambda_{u}^{p_{2}}$). 

If there is some element $\lambda\in \Lambda$ which cancels all $H_{1}(L;u)$,  then the element $(\lambda, -\lambda\mu(1-x))$ must be a boundary, and in particular $\lambda\in \mbox{Im}(\delta_{2}^{(1)})$. Therefore $\delta_{1}^{(1)}(\lambda)=0$, in other words $(1-x)\lambda =0$ in $\Lambda_{u}$. Now $(1-x)$ is invertible (as above), so $\lambda=0$ in $\Lambda_{u}$ and also in $\Lambda$, since $\Lambda\ri\Lambda_{u}$ is obviously injective. This ends the proof. 

\hfill $ \Box$

\vspace{.25in}

In conclusion, if $L$ is a monotone Lagrangian embedding which is displaceable, the condition $(a)$ implies that $H_{\ast}(L;u)$ is $\Lambda$-torsion for all classes $u$, whereas the condition $(b)$ implies that for some class it is not.  Therefore $L$ does not admit monotone Lagrangian embeddings which are displaceable. The first proof of \ref{general} is finished. 

\hfill $\Box$

\noindent\underline{Second proof}

  As in the first proof, we know using the hypothesis (a) that the lifted Floer homology is not well defined and therefore $\#_{2}^{g}\neq 0$ for some $g\in \piu(L)$. We will use the following 

\begin{propo}\label{conjugat} Let $L\subset M$ be an orientable monotone Lagrangian with Maslov number $N_{L}=2$. Suppose that for some $x\in L$ and some nontrivial $g\in \piu(L,x)$ we have $\#_{2}^{g}(x)\neq 0$. Then $g$ has a finite number of conjugacy classes in $\piu(L)$. In particular the subgroup $Z(g) \, =\, \{\gamma\in \piu(L,x)\, |\, \gamma g = g\gamma\}$ has finite index. 
\end{propo}

\noindent\underline{Proof} 

As pointed out in Section 1.2, given a path $\gamma$ between two $x,y\in L$, it yields a cobordism between $\mbox{ev}^{-1}(x)$ and $\mbox{ev}^{-1}(y)$ and if two disks passing through $x$ with boundaries in $g\in \piu(L,x)$  are respectively cobordant with two disks passing through $y$, the boundaries of these two disks define the same class in $\piu(L,y)$. Since $\#_{2}^{g}(x)=1$, there is at least one disk in $\mbox{ev}^{-1}(x)$ with boundary in $g$ which is cobordant with some disk in $\mbox{ev}^{-1}(y)$. 

Suppose now $x=y$. It is easy to see that the boundary of the latter disk is homotopic to $\gamma g \gamma^{-1}$ (homotopy with basepoint $x$). In particular, for every $\gamma\in \piu(L,x)$  there is some holomorphic disk in $\mbox{ev}^{-1}(x)$ with boundary in $\gamma g \gamma^{-1}$. On the other hand $\mbox{ev}^{-1}(x)$ has a finite number of elements, therefore $g$ has only a finite number of conjugacy classes. 

\hfill $\Box$

The end of the second proof of \ref{general} is obvious. If $\piu(L)$ is like in hypothesis $(b)$ then any non-trivial $g\in \piu(L)$ has an infinite number of conjugacy classes, contradicting the proposition above. Therefore $L$ cannot be monotone.

\hfill $\Box$

\vspace{.1in}

\noindent{\bf Remark} \\
The  statement \ref{conjugat} was inspired by a discussion with Jean-Fran\c{c}ois Barraud to whom I am grateful. 
\vspace{.25in}

\noindent\underline{Proof of \ref{connected}}

 Remark that $H_{n}(\covl_{j}, \bfz/2)=0$ implies that $\widetilde{L}_{j}$ are non-compact. A simple computation using Mayer-Vietoris shows that if $L_{1}$ and $L_{2}$ satisfy the hypothesis $(a)$ of Th.\ref{general} then so does $L=L_{1}\#L_{2}$. If $\beta_{1}(L_{j})\neq 0$ the proof follows directly from \ref{general} but actually we do not need this hypothesis. 

  According to Th.\ref{dadoi}, the lifted Floer homology $HF(\covl)$ is not well defined. By Prop.\ref{torsion} we get that all the Novikov homologies of $L$ are $\Lambda$-torsion. In particular this is the case for the usual homology $H(\covl,\bfz/2)$.  On the other hand, the computation of this homology shows that there is a free direct summand  $\Lambda$ in $H_{n-1}(\covl)$: a lift to $\covl$ of the sphere $\bfs^{n-1}$ along which the connected sum is done yields a generator of this summand. Therefore $H(\covl, \bfz/2)$ cannot be $\Lambda$-torsion and the proof is done.

\hfill $\Box$

\vspace{.25in}

\noindent\underline{Proof of \ref{product}}

By Gromov's theorem $H^{1}(L, \real)\neq 0$, so $L$ has an infinite fundamental group. Therefore $L$ satisfies the hypothesis $(a)$ of \ref{general} and as above we infer that $H(\covl, \bfz/2)$ is $\Lambda$-torsion. But this is impossible if $L$ is not prime, according to the previous proof. It follows that $L$ is prime and we can apply Fukaya's theorem Th.\ref{Fu} to finish our proof. 

\hfill $\Box$

\vspace{.1in} 

Let us now give a proof of \ref{product} which does not use Fukaya's theorem Th.\ref{Fu} which is technically difficult:\\

\noindent\underline{Second proof of \ref{product}}

As above, $L$ satisfies the hypothesis $(a)$, so lifted Floer homology is not well defined, which implies that there is some nontrivial $g\in \piu(L,x)$ such that the number of $2$-disks with boundary in $L$ passing through $x$ and belonging to $g$ is odd. Denote by $E=\{g_{1}, g_{2}, \ldots, g_{r}\}$ all the elements of $\piu(L,x)$ with this property. They satisfy the following conditions: 
\begin{itemize} 
\item (A) If $g\in E$ and $\gamma\in\piu(L,x) $ then $\gamma^{-1}g\gamma \in E$. In particular the subgroup $Z(g)$ of elements commuting with $g$ is of finite index in $\piu(L,x)$.
\item (B)  There is some morphism $u:\piu(L)\ri\bfz$ such that $u(g)=1$ for any $g\in E$.
\end{itemize}

  Indeed (A) is given by Prop.\ref{conjugat} and the morphism of (B) is $u = [\lambda]/2c$, where $\lambda$ is a primitive of the symplectic form on $\bfc^{3}$  and $c$ is the monotonicity constant of the manifold $L$. 

Let us prove first that $L$ fibers over the circle. According to a celebrated theorem of Stallings \cite{Sta} (for three dimensional manifolds) it is sufficient to prove that there is some morphism $\piu(L)\ri \bfz$ with finite generated kernel.

\begin{lem}\label{fib} \mbox{Ker}(u) is finitely generated. Therefore $L$ fibers over the circle. 
\end{lem}

\noindent\underline{Proof}

Let $x_{1}, x_{2}, \ldots, x_{s}$ a system of generators of $\piu(L)$. Denote $u(x_{i})=n_{i}\in \bfz$. Fix $t=g_{1}\in E$. We claim that $\mbox{Ker}(u)$ is spanned by $x_{i}t^{-n_{i}}$, $x_{i}^{-1}t^{n_{i}}$ and by $g_{j}g_{k}^{-1}$. Let $y\in \mbox{Ker}(u)$. \\
\\
1. $y\in H$, the subgroup spanned by $E$. 

We proceed by induction on $N(y)$, the number of letters $g_{i}^{\pm 1}$ in the writing of $y$. If $N=2$ then $y=g_{j}g_{k}^{-1}$ ($N=1$ is impossible). 

For a  general $y\in \mbox{Ker}(u)\cap H$ we have $y=zg_{j}g_{k}^{-1}v$ for some $z, v\in H$. According to property A, we have $zg_{j}g_{k}^{-1}=g_{j}g_{k}^{-1}z'$, with $z'\in H$ and $N(z')=N(z)$. The induction hypothesis applies to $z'v$. \\
\\
2. General case.  We have $y = x_{i}^{\pm 1}w=x_{i}^{\pm 1}t^{\bar{+} n_{i}}t^{\pm n_{i}}w$, for some $i$. Then, using the property A, $t^{\pm n_{i}}w=wh$, for some $h\in H$. Then we proceed in the same way with $w$  whose writing in the letters $x_{i}^{\pm 1}$ has one letter less than the one of $y$. Finally $y=y_{1}y_{2} $ where $y_{1}$  is a product of elements of the form $x_{i}t^{-n_{i}}$ and $x_{i}^{-1}t^{n_{i}}$  and $y_{2}\in \mbox{Ker}(u)\cap H$. We apply case 1 to $y_{2}$ and finish the proof.

\hfill $\Box$

\vspace{.1in} 

To complete the proof of \ref{product} we need:

\begin{lem}\label{pro}  Let $t=g_{1}\in E$ as above. Then $\piu(L)=Z(t)=\bfz \times \mbox{Ker}(u)$.
\end{lem}

\noindent\underline{Proof}

Consider the fibration given by the previous lemma. If the fiber $\Sigma$ is $\bfs^{2}$ there is nothing to prove. If not, suppose there is some $g\neq t$ in $E$. Therefore $Z(g)\cap Z(t)$ is of finite index in $\piu(L)$. By the property (B), $h= gt^{-1}\in \mbox{Ker}(u)=\piu(\Sigma)$. Obviously, $Z(h)\supset Z(g)\cap Z(t)$ is of finite index in $\piu(L)$ and therefore $Z(h)\cap \piu(\Sigma)$ is of finite index in $\piu(\Sigma)$. This is impossible unless $\Sigma$ is a torus. This means that when $\Sigma$ is of genus $\geq 2$ the set $E$ has only one element and $Z(t)=\piu(L)$. 

Suppose now $\Sigma =\bft^{2}$. Let $h_{1}$ and $h_{2}$ be generators of $\piu(\Sigma)$. There is some positive integer $m$ such that $h_{1}^{m}t=th_{1}^{m}$. On the other hand, we have $h_{1}t=th$ for some $h\in \mbox{Ker}(u)= \piu(\Sigma)$.  This implies $h_{1}^{m}=h^{m}$, therefore $h=h_{1}$ and $t$ commutes with $h_{1}$. Analogously $t$ commutes with $h_{2}$, so, again, $Z(t) = \piu(L)$. 

The application $\bfz\times \piu(\Sigma)\ri Z(t)=\piu(L)$, $(n,h)\mapsto t^{n}h$ is an isomorphism.

\hfill $\Box$

\vspace{.1in} 

We finish now the proof of \ref{product}. If $L\neq \bfs^{1}\times \bfs^{2}$ then $L$ is aspherical and  the genus of $\Sigma$ is greater then one. The previous lemma shows that the diffeomorphism  $\Gamma: \Sigma\ri \Sigma $ defined by the monodromy of the fibration over $\bfs^{1}$ induces the identity at the level of $\piu(\Sigma)$. The classical Dehn-Nielsen-Baer theorem \cite{Ni}, \cite{B} \cite{D}  implies then that $\Gamma$ is isotopic to the identity ant therefore  that $L$ is diffeomorphic to $\bfs^{1}\times \Sigma$. 

\hfill $\Box$

\vspace{.25in}

\noindent\underline{Proof of \ref{moregeneral}}

The point (i)  follows directly from Th.\ref{dadoi} and Prop.\ref{torsion}.\\
The point (ii) follows immediately from Th.\ref{dadoi} and Prop.\ref{conjugat}.

\hfill $\Box$

\vspace{.25in}

\noindent\underline{Proof of \ref{pair}}

Since $L_{1}$ and $L_{2}$ satisfy the condition $(a)$ of Th.\ref{general}, $L_{1}\times L_{2}$ has the same property, by K\"{u}nneth formula. Therefore according to \ref{dadoi}  the homology
$HF(\covl)$ is not defined, so using \ref{torsion} we get that $H(\covl, \bfz/2)$ is $\Lambda$-torsion. Again by K\"{u}nneth formula $H^{\mbox{dim}(L_{1})-1}(\covl_{1})\otimes_{\bfz/2}H^{\mbox{dim}(L_{2})-1}(\covl_{2})$ is a direct summand in $H^{\mbox{dim}(L)-2}(\covl)$. As in the proof of \ref{connected} there are direct summands $\bfz/2\, [\piu(L_{i})]$ in the homology groups of $\covl_{i}$ considered above. Their tensorial product over $\bfz/2$  $\Lambda=\bfz[\piu(L)]$ yields a direct summand in $H^{\mbox{dim}(L)-2}(\covl)$, therefore this homology cannot be $\Lambda$-torsion. This finishes the proof. 

\hfill $\Box$

\vspace{.25in}

\noindent\underline{Proof of \ref{fibration}}

Recall the formula (1) for the square of the lifted differential:  
$$(1)\gol (\delta^{\widetilde{L}})^{2} x\, =\, \sum_{g\in \piu(L)}\#_{2}^{g}\, (1+g)x. $$ 
The condition $\#_{2}=1$ implies that there is an odd number of terms in this sum. Denote by  $\alpha$  a primitive of the symplectic form (which was supposed to be exact) and let  $u=[\alpha|_{L}]$. We have by Stokes and the monotonicity of the submanifold $L$
$$u(g)\, =\, 2c$$ for any $g\in \piu(L)$ in the above formula, where $c>0$ is the monotonicity constant.  This implies
$$ (\delta^{\widetilde{L}})^{2} x\, =\, \lambda x\, = \,  (1+\mu)x,$$
where $u(g)>0$ for any $g$ appearing in the writing of $\mu$. Therefore $\lambda$ is invertible in $\Lambda_{u}$. On the other hand, we showed in the proof of \ref{torsion} that $\lambda$ cancels $H_{\ast}(L;u)$. It follows that this Novikov homology vanishes. 

\hfill $\Box$
 
\vspace{.25in}

\noindent\underline{Proof of \ref{spin}}

In \S 1.2 we defined Novikov homology with $\bfz/2$ coefficients. An analogous formula defines its version with integer coefficients $H_{\ast}(L,\bfz;u)$. For $u$ defined as in the previous proof, let us show that this Novikov homology is zero. When we work with integer coefficients, the formula $(1)$ has to be adapted. In the writing of $(\delta^{\widetilde{L}})^{2} x$ we should count {\it algebraically}  $\pm gx$ for each holomorphic disk with Maslov index two and boundary in $L$ of class $g$ passing through $x$ and $\pm x$ for each holomorphic disk with Maslov index two and boundary in $\phi_{1}(L)$ passing through $x$. The sum of the contributions of the disks of the second category is $\#=1$. Therefore, as above 
$$ (\delta^{\widetilde{L}})^{2} x\, = \,  (1+\mu)x,$$
where $u(g)>0$ for any $g$ appearing in the writing of $\mu\in \bfz[\piu(L)]$. As above, $1+\mu$ is invertible in $\Lambda_{u}$ and cancels $H_{\ast}(L,\bfz;u)$, so the latter vanishes.

The vanishing of $H_{\ast}(L,\bfz;u)$ is obviously a necessary condition for the existence of a fibration $f:L\ri \bfs^{1}$ in the class $u$ (actually in $u'=u/2c \in H^{1}(L,\bfz)$). In dimension greater than six there are two more conditions which together with $H_{\ast}(L,\bfz;u)=0$ are sufficient to the existence of this fibration \cite{La}, \cite{Pa}. The first one concerns the Whitehead torsion: $\tau(L;u)=0$ and the second one is $\mbox{Ker}(u)$ finitely presented. We remarked in \cite{Da1} (after the statement of Cor.3.6) that for $q$ odd we always have $\tau(L\times \bfs^{q};u)=0$. And if $\piu(L)$ si polycyclic, so is $\mbox{Ker}(u)$ and therefore it is finitely presented. This proves (i).

In order to prove (ii) choose $u'\in H^{1}(N,\bfz)\setminus\{0\}$ and consider the class $v=(u,u')$ on $L\times N$. We have $H_{\ast}(L\times N ; v)=0$. Indeed the identity of the Novikov complex $C_{\bullet}(L,u)$ is homotopic to zero (since this complex is free acyclic), and therefore the same property is valid for $C_{\bullet}(L\times N;v)=C_{\bullet}(L;u)\otimes_{\bfz}C_{\bullet}(N;u')$. Then, as above the torsion $\tau(L\times N\times \bfs^{3};v)$ is zero. Take a generic  closed one form on $L\times \bfs^{3}$ in the class $u$ with no zeros of Morse index $0$, $1$, $m-1$, $m$, where $m$ is the dimension of this manifold. This is always possible when $m\geq 5$ and the Novikov homology vanishes, as pointed out by F. Latour \cite{La}. Then take a generic closed one-form on $N$  belonging to $u'$ without minimum and maximum (which is always possible in a non zero cohomology class). Their sum has no zero of Morse index and co-index $\leq 2$. This enables one to cancel all its zeros as in \cite{La} and to prove that it is cohomologous with a non singular one. Therefore $L\times N\times \bfs^{3}$ fibers over the circle and the proof is finished.

\hfill $\Box$ 

\vspace{.25in}

\noindent {\bf Acknowledgements} I thank Jean-Fran\c{c}ois Barraud, Gw\'ena\"{e}l Massuyeau and Jean-Claude Sikorav for our valuable discussions on the subject.

\end{document}